\documentclass[10pt]{article}
\usepackage{amsfonts}
\usepackage{mathrsfs}
\usepackage{amsmath}
\usepackage{amssymb}
\usepackage{graphicx}
\graphicspath{{figure/}}
\usepackage{epic}
\usepackage{amsthm,latexsym}
\renewcommand{\paragraph}{\roman{paragraph}}
\def \G{\Gamma}
\DeclareMathOperator{\sgn}{sgn}
\newtheorem{theorem}{\scshape \mdseries  Theorem}[section]
\newtheorem{lemma}[theorem]{\scshape \mdseries  Lemma}
\newtheorem{coro}[theorem]{\scshape \mdseries  Corollary}

\begin{document}

\title{The Nullity of  Bicyclic Signed Graphs\thanks{
      Supported by National Natural Science Foundation of China (11071002, 11126178), Program for New Century Excellent
Talents in University, Key Project of Chinese Ministry of Education
(210091), Specialized Research Fund for the Doctoral Program of
Higher Education (20103401110002), Science and Technological Fund of
Anhui Province for Outstanding Youth (10040606Y33), Project of Educational Department of Anhui Province (KJ2011A019),
Scientific
Research Fund for Fostering Distinguished Young Scholars of Anhui
University(KJJQ1001), Academic Innovation Team of Anhui University
Project (KJTD001B).}}
\author{Yi-Zheng Fan\thanks{Corresponding author. E-mail addresses: fanyz@ahu.edu.cn (Y.-Z. Fan), wenxuedu@gmail.com (W.-X. Du), 786189560@qq.com (C.-L. Dong).}, \
Wen-Xue Du, \  Chun-Long Dong\\
 {\small  \it School of Mathematical Sciences, Anhui University, Hefei 230601, P. R. China}
 }
\date{}
\maketitle

\noindent{\bf Abstract:} Let $\G$ be a signed graph and let $A(\G)$ be the adjacency matrix of $\G$.
The nullity of $\G$ is the multiplicity of eigenvalue zero in the spectrum of $A(\G)$.
In this paper we characterize the signed graphs of order $n$ with nullity $n-2$ or $n-3$, and introduce a graph transformation which preserves the nullity.
As an application we determine the unbalanced bicyclic signed graphs of order $n$ with nullity $n-3$ or $n-4$,
and signed bicyclic signed graphs (including simple bicyclic graphs) of order $n$ with nullity $n-5$.


\noindent{\bf MR(2000) Subject Classification:} 05C50

\section{Introduction}
A {\it signed graph} is a graph with a sign attached to each of its edges.
Formally, a signed graph $\Gamma=(G,\sigma)$ consists of a simple graph $G=(V,E)$,
  referred to as its underlying graph, and a mapping $\sigma: E \rightarrow \{+,-\}$, the edge labeling.
To avoid confusion, we also write $V(\G)$ or $V(G)$ instead of $V$, $E(G)$ instead of $E$, and $E(\G)=E^\sigma$.
The underlying graph $G$ of $\G$ is sometimes denoted by $^u\G$.

The {\it adjacency matrix} of  $\G$ is  $A(\Gamma)=(a_{ij}^\sigma)$ with
$a_{ij}^\sigma=\sigma(v_iv_j)a_{ij}$, where $(a_{ij})$ is the adjacency matrix of the underlying graph $G$.
In the case of $\sigma=+$, which is an all-positive edge labeling,
   $A(G,+)$ is exactly the classical adjacency matrix of $G$.
So a {\it simple graph} is always assumed as a signed graphs with all edges positive.
The {\it nullity of a signed graph} $\G$ is defined as the multiplicity of the eigenvalue zero in the spectrum of $A(\G)$,
  and is denoted  by $\eta(\G)$.
The {\it rank} of $\G$ is referred to the rank of $A(\G)$, and is denoted by $r(\G)$. Surely $\eta(\G)+r(\G)=n$ if $\G$ has $n$ vertices.

 Let $C$ be a cycle of $\Gamma$.
The {\it sign} of $C$, denoted by $\sgn(C)$, and is defined as $\sgn(C)=\Pi_{e\in C} \sigma(e)$.
The cycle $C$ is said to be {\it positive} or {\it negative} if $\sgn(C)=+$ or $\sgn(C)=-$.
A signed graph is said to be {\it balanced} if all its cycles are positive,
  or equivalently, all cycles have an even number of negative edges; otherwise it is called {\it unbalanced}.

  Suppose $\theta: V(G) \to \{+,-\}$ is any sign function.
Switching $\G$ by $\theta$ means forming a new signed graph $\G^\theta=(G, \sigma^\theta)$ whose underlying
graph is the same as $G$, but whose sign function is defined on an edge $uv$ by
$\sigma^\theta(uv)=\theta(u)\sigma(uv)\theta(v)$.
Note that switching does not change the signs or balanceness of the cycles of $\G$.
If we define a (diagonal) signature matrix $D^\theta$ with $d_v=\theta(v)$ for each $v \in V(G)$, then
$A(\G^\theta)=D^\theta A(\G)D^\theta$.
Two graphs $\G_1,\G_2$ are called {\it switching equivalent}, denoted by $\G_1 \sim \G_2$, if there exists a switching function $\theta$ such that
  $\G_2=\G_1^\theta$, or equivalently $A(\G_2)=D^\theta A(\G_1)D^\theta$.

\begin{theorem} {\em \cite{hou}} \label{balance}
Let $\G$ be a signed graph.
Then $\Gamma$ is balanced if and only if $\G=(G, \sigma) \sim (G,+)$.
\end{theorem}

Recently the authors \cite{fanww} investigate the nullity of signed graphs, and determine the unicyclic signed graphs with nullity $n-2$, $n-3$, $n-4$, or $n-5$.
In this paper we characterize the signed graphs of order $n$ with nullity $n-2$ or $n-3$, and introduce a graph transformation which preserves the nullity.
As an application we determine the unbalanced bicyclic signed graphs of order $n$ with nullity $n-3$ or $n-4$,
and signed bicyclic signed graphs (including simple bicyclic graphs) of order $n$ with nullity $n-5$.
We hope this work would be helpful to characterize the structure of signed graphs from the view of nullity.

Finally we note some work on the nullity of simple graphs.
For a simple graph $G$ of order $n$ containing at least one edge, $0\le \eta(G)\le n-2$.
Cheng and Liu \cite{chengl} characterize the simple graphs of order $n$ with nullity $n-2$ or $n-3$.
Cheng, Huang and Yeh \cite{chenghy,chenghy2} characterize the simple graphs of order $n$ with nullity $n-4$ or $n-5$.
The characterization of simple graphs of order $n$ with nullity $n-6$ or more is still open.
Much work is devoted to the nullity of special classes of simple graphs; see \cite{fan, fiorini, gong, gongxu, guo, gut, hu, lics, li, lic, nath, tan, zhu}.
The work on the nullity of simple graphs were found to be useful in chemistry to investigate the stability of the alternant hydrocarbons (\cite{collatz},\cite{long}),
and is also of interest in mathematics itself, as it is closely related to the minimum rank problem of symmetric
matrices whose patterns are described by graphs \cite{fallat}.

\section{Nullity of signed graphs}
Some notations are listed below.
Let $\G$ be a signed graph.
Let $v$ be a vertex of $\G$ and let $U \subset V(\G)$.
Denote by $N(v)$ the neighborhood of $v$ in $\G$, $d(v)$ the degree of $v$ (i.e. the cardinality of $N(v)$).
Denote by $\G[U]$ the induced subgraph of $\G$ on the vertices of $U$ including the signs of edges.
Denote by $P_r,C_r,K_r$ respectively a path, a cycle and a complete graph, all of which are simple graphs of order $r$.
We use the notation $G_1 + G_2$ to denote the union of two vertex-disjoint simple graphs $G_1,G_2$.

\begin{lemma} {\em \cite{fanww}} \label{tree-cycle}
{\em (i)}  If $T$ is an acyclic signed graph or a signed tree of order $n$. Then
$\eta(T)=n-2\mu(T)$.

{\em (ii)}  If $C_n$ is a balanced cycle,
then $\eta(C_n)=2$ if $n \equiv 0(\!\!\!\mod 4)$ and $\eta(C_n) = 0$ otherwise.

{\em (iii)}  If $C_n$ is an unbalanced signed cycle,
then $\eta(C_n)=2$ if $n \equiv 2 (\!\!\!\mod 4)$, and $\eta(C_n)=0$ otherwise.
\end{lemma}

If a signed graph $\G$ contains at least one edge, then $\eta(\G) \le n-2$.
By adopting the technique in \cite{chengl}, We first characterize the signed graphs of order $n$ with nullity $n-2$ or $n-3$.

\begin{lemma} \label{rank3}
Suppose that $\G$ is a signed graph of order $n$ and $\G$ has no isolated vertex.
Let $x$ be an arbitrary vertex in $\G$. Let $Y=N(x)$ and $X = V(\G)-Y$.
If $r(\G) \le 3$, then

{\em (i)} No two vertices in $X$ are adjacent.

{\em (ii)} Each vertex from $X$ and each vertex from $Y$ are adjacent.

\end{lemma}

\noindent {\bf Proof.}
Let $x_1,x_2$ be two vertices in $X$.
If one of them is $x$, surely the assertion (i) holds.
Now assume $x_1$ and $x_2$ are adjacent and neither of them is $x$.
Since $\G$ contains no isolated vertices, $Y$ is not empty.
Select any vertex $y$ in $Y$. Then $^u\G[x_1,x_2,y]$ is isomorphic to $P_2 + K_1$, $P_3$ or $K_3$,
and hence $^u\G[x,x_1,x_2,y]$ is isomorphic to $P_2+P_2$, $P_4$ or $K_3$ with a pendant edge.
In any case the subgraph $\G[x,x_1,x_2,y]$ has rank $4$, a contradiction.

For the assertion (ii), assume to the contrary,
there exist $x_1 \in X, y_1 \in Y$ such that $x_1$ is not adjacent to $y_1$.
Surely $x_1 \ne x$.
As $\G$ has no isolated vertices, $x_1$ has a  neighbor $z$, necessarily in $Y$ by the result (i).
The subgraph $\G[x,x_1,y_1,z]$ has rank $4$ whether $y_1$ and $z$ are adjacent or not, also a contradiction.
\hfill$\blacksquare$

 \begin{theorem} \label{rank2}
Let $\G$ be a signed graph of order $n \ge 2$.
Then  $\eta(\G) = n-2$ or equivalently $r(\G)=2$ if and only if $\G$ is a balanced complete bipartite graph together with some isolated vertices.
 \end{theorem}

\noindent {\bf Proof.}
The sufficiency is clear. So we only prove the necessity.
First assume $\G$ has no isolated vertices.
Choose an arbitrary vertex $x$ in $G$.
Let $Y=N(x), X= V(\G)-Y$, both of which are nonempty.
We assert any two vertices of $Y$ are not adjacent; otherwise $\G$ contains a triangle which has rank $3$ by Lemma \ref{tree-cycle}.
Then $\eta(\G) \le n-3$, a contradiction.
Hence $\G$ is a complete bipartite graph by Lemma \ref{rank3}.

Next we show $\G$ is balanced.
Otherwise $\G$ contains an unbalanced cycle, say $C_{2k}$.
Further assume $C_{2k}$ has the minimum length.
By Lemma \ref{tree-cycle}(iii), $k \ge 3$, as an unbalanced $C_4$ would have rank $4$.
Adding an edge between two vertices on $C_{2k}$ with distance $3$,
the cycle is split into two cycles $C_4$ and $C_{2k-2}$, both being subgraphs of $\G$.
By the assumption on $C_{2k}$, $C_4$ is balanced, and hence $C_{2k-2}$ is unbalanced, a contradiction.

Finally, if $\G$ contains isolated vertices, letting $H$ be the subgraph of $\G$ by deleting those isolated vertices,
then $r(H)=r(\G)=2$.
So $H$ is a balanced complete bipartite graph.
The result now follows.\hfill$\blacksquare$

Let $G$ be a signed graph and let $v$ be a vertex of $G$.
The {\it positive neighborhood} (respectively, {\it negative neighborhood}) of $v$ is the set of vertices joining $v$
  with positive edges (respectively, negative edges).

 \begin{theorem} \label{rank-3}
Let $\G$ be a signed graph of order $n \ge 2$.
Then  $\eta(\G) = n-3$ or equivalently $r(\G)=3$ if and only if $\G$ is a complete tripartite graph together with some isolated vertices,
such that all vertices in each of three parts of the complete tripartite graph have the same positive neighborhoods and negative neighborhoods.
 \end{theorem}

\noindent {\bf Proof.}
If $\G$ is the graph as in the theorem, then $\G$ contains a triangle, and by Lemma \ref{tree-cycle}, $r(\G) \ge 3$.
From the matrix of $A(\G)$, all rows corresponding the vertices in a given part of the complete tripartite graph are same, which implies that $r(\G) \le 3$.
So the sufficiency follows.

Now we prove the necessity.
 First assume $\G$ has no isolated vertices.
Choose an arbitrary vertex $x$ in $G$.
Let $Y=N(x), X= V(\G)-Y$, both of which are nonempty.
By Lemma \ref{rank3}, any two vertices in $X$ are not adjacent, and any vertex from $X$ is adjacent to any vertex from $Y$.

We assert $r(\G-X)=2$.
If $r(\G-X)>2$, then $r(G)=r(\G-X)=3$.
So $\G-X$ contains an induced subgraph of order $3$ (i.e. a triangle) with rank $3$.
Therefore, $G$ contains $K_4$ as a subgraph. If $K_4$ contains an unbalanced $C_4$, by Lemma 2.1(iii), $r(K_4)=4$, a contradiction.
Otherwise, by a little calculation we also get $r(K_4)=4$, a contradiction.
So $r(\G-X)\le 2$.
If $r(\G-X) <2$, then $r(\G-X) =0$, and hence $\G-X$ contains no edges.
So, $\G$ is a complete bipartite graph, whose rank is an even number, a contradiction.

From above discussion, $\G-X$ is a balanced complete bipartite graph together with some isolated vertices by Theorem \ref{rank2}.
If $\G-X$ contains isolated vertices, then $\G$ contains a triangle with an pendant edge as an induced subgraph.
By Lemma \ref{pend} below, the latter subgraph has rank $4$, a contradiction.
So $\G$ is a complete tripartite graph.
Write $A(\G)$ as the following form:
$$A(\G)=\left[\begin{array}{ccc}
O & J_{12} & J_{13} \\
J_{12}^T & O & J_{23} \\
J_{13}^T & J_{23}^T & O
\end{array}\right]=:\left[\begin{array}{c}
A_1 \\
A_2 \\
A_3\end{array}\right],
$$
where each $J_{ij}\;(i \ne j)$ consists of entries $1$ or $-1$.
If there exists a row submatrix $A_i$, say $A_1$, has rank $r \ge 2$, then $A_1^T$ also has rank $r$, which implies that $r(\G) \ge 2r \ge 4$, a  contradiction.
So, each row of $A_i$ is a copy of the first row of $A_i$ for $i=1,2,3$. The result follows in this case.

If $\G$ contains isolated vertices, letting $H$ be the subgraph of $\G$ by deleting those isolated vertices,
then $r(H)=r(\G)=3$.
So $H$ has the property as discussed above. The result follows.\hfill$\blacksquare$

At the end of this section, we introduce some graph transformations which preserve the nullity.

\begin{lemma} {\em \cite{fanww}} \label{pend}
Let $\G$ be a signed graph containing a pendant vertex, and let $H$ be the induced subgraph
of $G$ obtained by deleting this pendant vertex together with the vertex adjacent to it. Then
$$\eta(\G)=\eta(H).$$
\end{lemma}

\begin{coro} \label{pendc} Let $\G$ be a signed graph of order $n \ge 4$ which is not a star.
If $\G$ contains a pendant vertex, then $\eta(\G)\leq n-4$.
\end{coro}

\noindent{\bf Proof.}
Let $v$ be a pendant vertex of $\G$, and let $u$ be the vertex adjacent to $v$.
By Lemma \ref{pend}, $\eta(\G)=\eta(\G-v-u)$, and hence $r(\G)=r(\G-u-v)+2$.
As $\G$ is not a star, $\G-u-v$ contains at least one edge so that $ r(\G-u-v) \ge 2$, which implies $r(\G) \ge 4$, or equivalently $\eta(\G)\leq n-4$.
\hfill$\blacksquare$

Denote by $Pv_1v_2v_3$ a path on vertices $v_1,v_2,v_3$ with edges $v_1v_2,v_2v_3$.
If a signed graph $\G$ contains a path $Pv_1v_2v_3$ which satisfies $d(v_2)=2, v_1v_3 \notin E(\G), N(v_1) \cap N(v_3) \setminus \{v_2\} = \emptyset$,
then $Pv_1v_2v_3$ is called a {\it special path} of $\G$.

\begin{lemma} \label{spath1} Let $\G$ be a signed graph containing a special path $Pv_{1}v_{2}v_{3}$, where $\sigma(v_{1}v_{2})=-1, \sigma(v_2v_3)=1$.
Suppose that $N(v_1)\setminus \{v_2\} \ne \emptyset$ and $v \in N(v_1)\setminus \{v_2\}$.
Let $\G'$ be obtained from $\G$ by deleting the edge $vv_1$ and adding a new edge $vv_3$ with the same sign of $vv_1$.
Then $\eta(\G')=\eta(\G)$.
\end{lemma}

\noindent{\bf Proof.}
Without loss of generality, write $A(\G)$ as follows, where the first to the fourth rows correspond to the vertices $v_1,v_2,v_3,v$ respectively,
and $H=\G-\{v_1,v_2,v_3,v\}$.
$$A(\G)=\left[
         \begin{array}{ccccc}
           0 & -1 & 0 & \sigma(vv_1) & \alpha\\
           -1 & 0 & 1 & 0 & \textbf{0}\\
           0 & 1 & 0 & 0 & \beta\\
           \sigma(vv_1) & 0 & 0 & 0 & \gamma\\
           \alpha^{T} & \textbf{0}^T & \beta^{T} & \gamma^{T} & A(H)\\
         \end{array}
       \right]
$$
Taking the second row (column) to act on the fourth row (column) by elementary transformations, we arrive at a matrix:
$$ M=\left[
         \begin{array}{ccccc}
           0 & -1 & 0 & 0 & \alpha \\
           -1 & 0 & 1 & 0 & \textbf{0} \\
           0 & 1 & 0 & \sigma(vv_1) & \beta \\
           0 & 0 & \sigma(vv_1) & 0 & \gamma \\
           \alpha^{T} & \textbf{0}^T & \beta^{T} & \gamma^{T} & A(H) \\
         \end{array}
       \right]
$$
Evidently $M$ has the same rank as $A(\G)$, is exactly the adjacency matrix of $\G'$. The result follows. \hfill$\blacksquare$

\begin{coro} \label{spath2} Let $\G$ be a signed graph containing a special path $Pv_{1}v_{2}v_{3}$, where $\sigma(v_{1}v_{2})=-1, \sigma(v_2v_3)=1$.
 Let $\G'$ be a signed graph obtained from $\G$ by contracting the path into a single vertex.
Then $\eta(\G')=\eta(\G)$.
\end{coro}

\noindent{\bf Proof.}
If $N(v_1)\setminus \{v_2\} \ne \emptyset$, deleting the edge $vv_1$ and adding a new edge $vv_3$ with same signs as $vv_1$ for each vertex $v \in N(v_1)\setminus \{v_2\}$,
we will arrive at a graph $H$ such that $\eta(H)=\eta(\G)$ by Lemma \ref{spath1}.
Noting that $v_1$ is a pendant vertex of $H$, by Lemma \ref{pend}, $\eta(H)=\eta(H-v_1-v_2)$.
It is readily seen that $H-v_1-v_2$ is isomorphic to $\G'$. The result follow.\hfill$\blacksquare$

A {\it bicyclic graph} $G$ is a simple connected graph such that $|E(G)|=|V(G)|+1$.

\begin{coro} \label{spathbic} Let $B$ be a bicyclic signed graph of order $n$.
If $B$ contains a special path, then $\eta(B)\leq n-4$.
\end{coro}

\noindent{\bf Proof.}
Assume $Pv_{1}v_{2}v_{3}$ is a special path of $B$.
There exists a sign function $\theta$ such that in the graph $B^\theta$ the sign of $v_1v_2$ is negative and the sign of $v_2v_3$ is positive.
Noting that $\eta(B)=\eta(B^\theta)$, so we deal with the graph $B^\theta$.
By Corollary \ref{spath2}, we can get a new signed graph $B'$ from $B^\theta$ such that $\eta(B')=\eta(B^\theta)=\eta(B)$ and $r(B)=r(B^\theta)=r(B')+2$, where $B'$ is obtained from $B^\theta$ by contracting  the special path into a single vertex.
As $B$ is bicyclic, $B'$ contains at least one edge, $r(B') \ge 2$ and hence $r(B) \ge 4$.
The result follows. \hfill$\blacksquare$

\section{Nullity of bicyclic signed graphs}

There are two basic bicyclic graphs: $\infty$-graph and $\Theta$-graph.
An {\it $\infty$-graph}, denoted by $\infty(p,q,l)$, is obtained from two vertex-disjoint cycles $C_p$ and $C_q$ by
connecting one vertex of $C_p$ and one of $C_q$ with a path $P_l$ of length $l-1$ (in the case of $l = 1$, identifying the above two vertices); and
a {\it $\Theta$-graph}, denoted by $\Theta(p,q,l)$, is a union of three internally disjoint paths $P_{p+1},P_{q+1},P_{l+1}$ of length $p,q,l$ respectively with common end vertices, where $p, q, l \ge 1$ and at most one of them is $1$.
Observe that any bicyclic graph $G$ is obtained from an $\infty$-graph or a $\Theta$-graph (possibly) by attaching trees to some of its vertices.

Note that the simple bicyclic graphs (or balanced signed bicyclic graphs) of order $n$ with nullity $n-3$ or $n-4$ have been characterized in \cite{hu}.
In this section we will characterize the unbalanced signed bicyclic graphs of order $n$ with nullity $n-3$ or $n-3$, 
and the signed bicyclic signed graphs (including simple bicyclic graphs) of order $n$ with nullity $n-5$..

\begin{theorem} \label{main1}
Let $B$ be an unbalanced bicyclic signed graph of order $n$.
Then $\eta(B) \le n-3$, with equality if and only if $^uB=\Theta(2,2,1)$ and the two triangles of $B$ are both unbalanced.
\end{theorem}

\noindent{\bf Proof.}
%
As $B$ is unbalanced, by Theorem \ref{rank2}, $\eta(B) \le n-3$.
If $\eta(B) = n-3$, by Theorem \ref{rank-3}, $B$ is a complete tripartite graph, which implies $^uB=\Theta(2,2,1)$, and two triangles of $B$ are both unbalanced.
The sufficiency also follows from Theorem \ref{rank-3}.\hfill$\blacksquare$


\begin{center}
\includegraphics[scale=0.5]{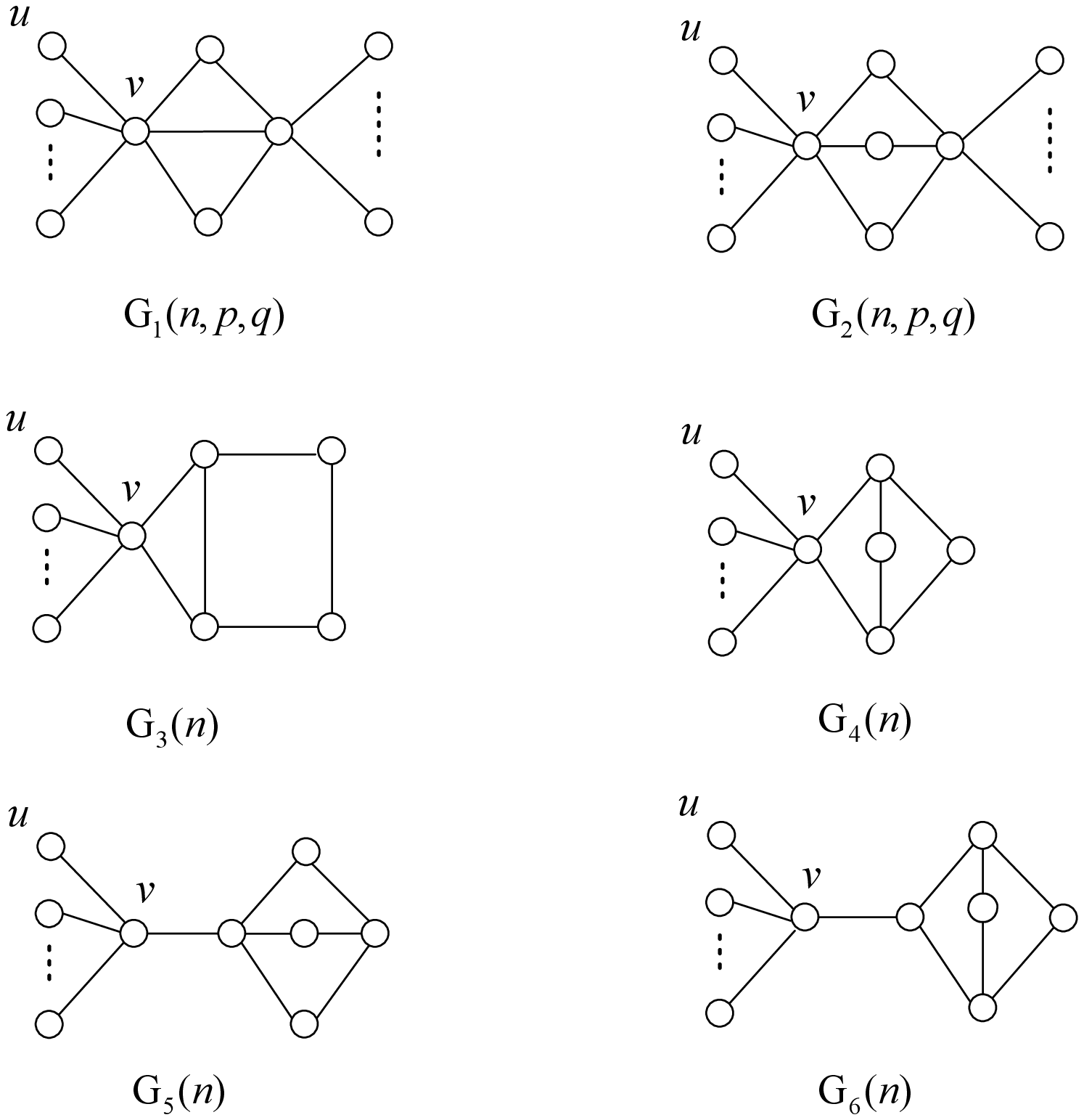}

{\small Fig. 3.1. Six bicyclic graphs with pendant vertices}
\end{center}

\begin{theorem} \label{main2}
Let $B$ be an unbalanced bicyclic signed graph of order $n$.
Then $\eta(B)= n-4$  if and only if $B$ is one of the following graphs with certain properties:

{\em (1)} $^uB=G_1(n,p,q)$ of Fig. 3.1, and $B$ contains at least one unbalanced cycle,

{\em (2)} $^uB=G_2(n,p,q)$ of Fig. 3.1, and $B$ contains at least one unbalanced cycle,

{\em (3)} $^uB=G_3(n)$ of Fig. 3.1,  and the triangle of $B$ is unbalanced and the quadrangle of $B$ is balanced,

{\em (4)} $^uB=G_4(n)$ of Fig. 3.1,  and the quadrangles of $B$ containing the vertex $v$ are unbalanced and the quadrangle of $B$ containing no $v$ is balanced,

{\em (5)} $^uB=G_{12}$ of Fig. 3.2, and $B$ contains exactly one unbalanced cycle,

{\em (6)} $^uB=G_{14}$ of Fig. 3.3, and $B$ contains an unbalanced triangle and a balanced quadrangle,

{\em (7)} $^uB=G_{15}$ of Fig. 3.3, and $B$ has exactly one unbalanced triangle,

{\em (8)} $^uB=G_{16}$ of Fig. 3.3, and $B^\theta$ has an edge labeling as in Fig. 3.3 for some sign function $\theta$.
\end{theorem}

{\bf Proof.}
The sufficiency can be verified by Lemma \ref{pend}, Lemma \ref{tree-cycle}, Theorem \ref{rank2}, Corollary \ref{spath2} or direct calculations.

Now assume $\eta(B)= n-4$. We divide the discussion into cases.

Case 1: $B$ contains pendant vertices.
Let $u$ be a pendant vertex of $B$, and let $v$ be the vertex adjacent to $u$.
By Lemma \ref{pend}, $\eta(B)=\eta(B-u-v)=n-4$.
So $r(B-u-v)=2$.
By Theorem \ref{rank2}, $B-u-v$ is a balanced complete graph $H$ together with some isolated vertices.
As $B$ is bicyclic,  $^uH$ is a star, or $C_4$ or $\Theta(2,2,2)$.
If $^uH=\Theta(2,2,2)$, then $^uB$ is $G_5(n)$ or $G_6(n)$ of Fig. 3.1.
In this case $B$ is also balanced as $H$ is balanced. So this case cannot occur.
If $H$ is a star, then $^uB$ is $G_1(n,p,q)$ or $G_2(n,p,q)$ of Fig. 3.1, and $B$ contains at least one unbalanced cycle.
If $H=C_4$, then $^uB$ is $G_3(n)$ or $G_4(n)$  of Fig. 3.1.
In addition, if $^uB=G_3(n)$, then the triangle is unbalanced and the quadrangle is balanced;
if $^uB=G_4(n)$, then the  quadrangle containing no $v$ is balanced and the quadrangles containing $v$ are unbalanced.

Case 2: $B$ contains no pendant vertices.
First assume that $B$ is an $\infty$-graph and contains a special path $P$.
Without loss of generality, assume $P$ contains a positive edge and a negative edge.
Now contracting the path $P$ into a single vertex, by Corollary \ref{spath2},
the resulting graph, denoted by $B'$, is still an $\infty$-graph, has the same nullity as $B$, and hence has rank $2$.
By Theorem \ref{rank2}, $B'$ is a balanced complete bipartite graph,  which is impossible.
So $B$ contains no special paths, and $^uB$ is one of graphs $G_i, i=7,8,\ldots,12$, of Fig. 3.2.
If $^uB$ is one of $G_7,G_8,G_9$, as $B$ is unbalanced, $B$ contains $C_3+P_2$ or $C_4+P_2$ as an induced subgraph, where $C_3$ or $C_4$ is unbalanced.
By Lemma \ref{tree-cycle}, $r(C_3+P_2)=5$, $r(C_4+P_2)=6$, a contradiction.
So, it suffices to consider the case of $^uB$ being one of $G_{10},G_{11},G_{12}$.
By a direct calculation, $^uB=G_{12}$, and $B$ contains exactly one unbalanced cycle.

\begin{center}

\includegraphics[scale=0.5]{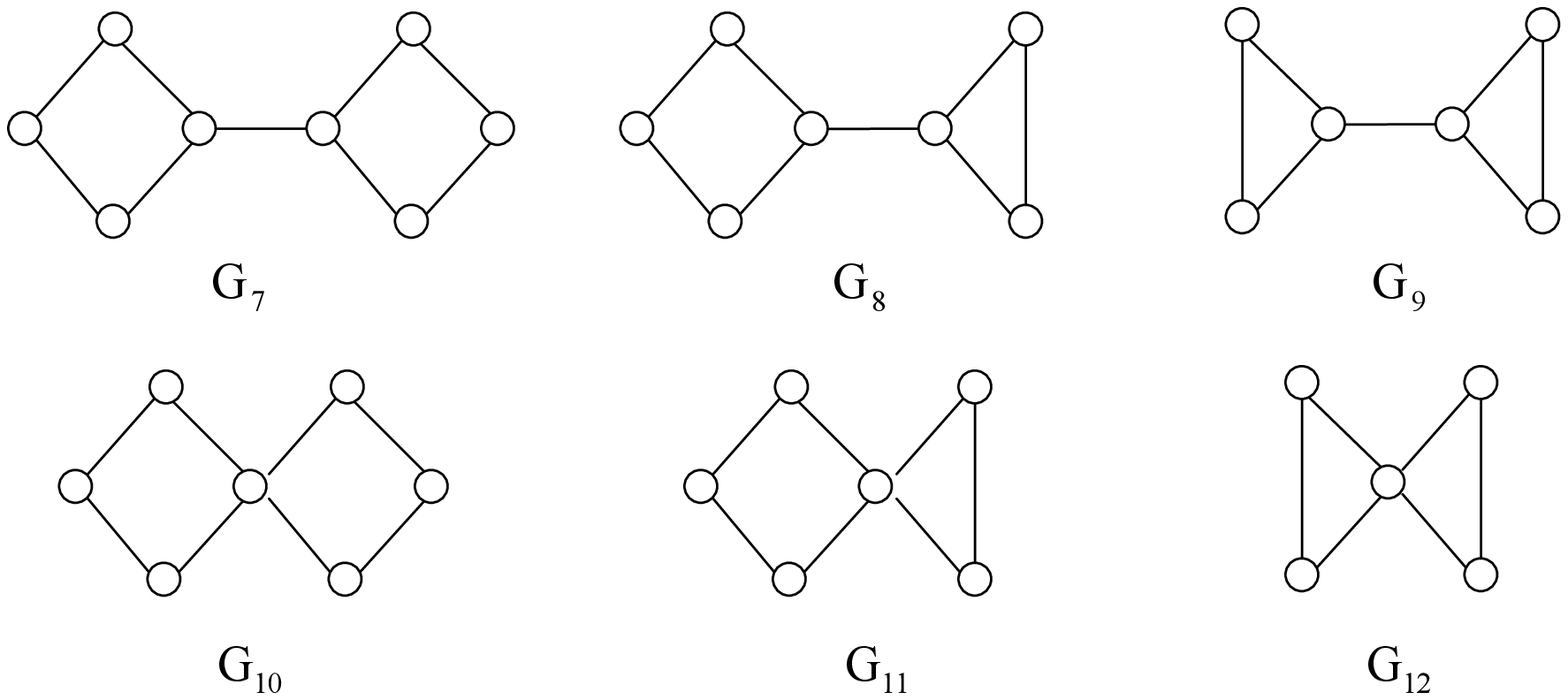}

{\small Fig. 3.2. Six $\infty$-graphs}

\end{center}

Next assume that $B$ is a $\Theta$-graph.
If $B$ contains no special paths, then $^uB$ is one of the graphs $G_{13},G_{14},G_{15}$ of Fig. 3.3.
By a direct calculation, $^uB=G_{14}$ and $B$ has an unbalanced triangle and a balanced quadrangle, or
$^uB=G_{15}$ and $B$ has exactly one unbalanced triangle.
If $B$ contains a special path $P$, without loss of generality, $P$ starts at a vertex of degree $3$, say $Pv_1v_2v_3$; see $G_{16}$ of Fig. 3.3.
Without loss of generality, the edge $v_1v_2$ is negative and $v_2v_3$ is positive.
Now contracting the path $P$ into a single vertex say $u$,  by Corollary \ref{spath2}, the resulting graph $\bar{B}$ will have rank $2$.
By Theorem \ref{rank2}, $\bar{B}$ is a balanced complete bipartite graph, that is $^u\bar{B}=\Theta(2,2,2)$.
By Theorem \ref{balance}, there exists a sign function $\eta$ defined on $V(\bar{B})$, such that all edges of $\bar{B}^\eta$ are positive.
Returning to the original graph, $^uB=G_{16}$ of Fig. 3.3.
Define a sign function $\theta$ on $B$ such that $\theta(v_1)=\theta(v_2)=\theta(v_3)=\eta(u)$, and $\theta(v)=\eta(v)$ for any other vertex $v$.
Then $B^\theta$ is the graph $G_{16}$ with edge labeling as in Fig. 3.3.\hfill$\blacksquare$

\begin{center}

\includegraphics[scale=0.5]{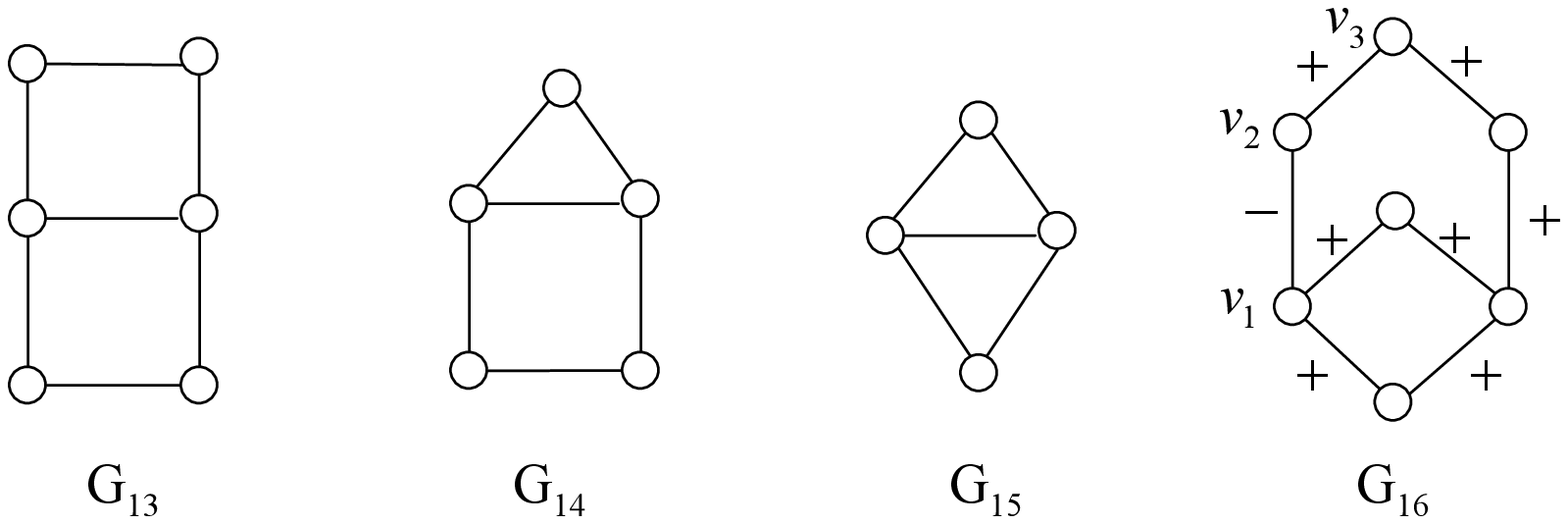}

{\small Fig. 3.3. Four $\Theta$-graphs}
\end{center}

\begin{theorem} \label{main3}
Let $B$ be a bicyclic signed graph of order $n$.
Then $\eta(B)= n-5$  if and only if $B$ is one of the following graphs with certain properties:

{\em (1)} $^uB=G_{17}$ of Fig. 3.4;

{\em (2)} $^uB=G_{18}$ of Fig. 3.4, and the two triangles of $B$ have the same balanceness;

{\em (3)} $^uB=G_{19}$ of Fig. 3.4, and the two triangles of $B$ have the same balanceness;

{\em (4)} $^uB=G_{11}$ of Fig. 3.2, and $B$ is balanced;

{\em (5)} $^uB=G_{12}$ of Fig. 3.2, and the two triangles of $B$ have the same balanceness;

{\em (6)} $^uB=G_{14}$ of Fig. 3.3, and the quadrangle of $B$ is unbalanced;

{\em (7)} $^uB=G_{20}$ or $G_{21}$ of Fig. 3.4, $B^\theta$ has an edge labeling as in Fig. 3.4 for some sign function $\theta$.

\end{theorem}

{\bf Proof.}
The sufficiency can be verified by Lemma \ref{pend}, Lemma \ref{tree-cycle}, Theorem \ref{rank-3}, Corollary \ref{spath2} or direct calculations.

Now assume that $\eta(B)=n-5$. We divide the discussion into cases.

Case 1: $B$ contains pendant vertices.
Let $u$ be a pendant vertex of $B$, and let $v$ be the vertex adjacent to $u$.
By Lemma \ref{pend}, $\eta(B)=\eta(B-u-v)=n-5$.
So $r(B-u-v)=3$.
By Theorem \ref{rank-3}, $B-u-v$ is a complete tripartite graph $H$ together with some isolated vertices.
As $B$ is bicyclic,  $^uH$ is a triangle, or $\Theta(2,2,1)$.
Hence $^uB$ is one of the graphs $G_{17},G_{18},G_{19}$ in Fig. 3.4, where the two triangles of $B$ have the same balanceness for the last two cases by Theorem \ref{rank-3}.

Case 2: $B$ contains no pendant vertices.
First assume that $B$ is an $\infty$-graph and contains a special path $P$.
Without loss of generality, assume $P$ contains a positive edge and a negative edge.
Now contracting the path $P$ into a single vertex, by Corollary \ref{spath2},
the resulting graph, denoted by $B'$, is still an $\infty$-graph, has the same nullity as $B$, and hence has rank $3$.
By Theorem \ref{rank2}, $B'$ is a complete tripartite graph,  which is impossible.
So $B$ contains no special paths, and $^uB$ is one of graphs $G_i, i=7,8,\ldots,12$, of Fig. 3.2.
By a little calculation, $^uB=G_{11}$ and $B$ is balanced, or  $^uB=G_{12}$ and two triangles of $B$ have the same balanceness.

Next assume that $B$ is a $\Theta$-graph.
If $B$ contains no special paths, then $^uB$ is one of the graphs $G_{13},G_{14}$ of Fig. 3.3.
By a direct calculation, $^uB=G_{14}$ and the quadrangle of $B$ is unbalanced.
If $B$ contains a special path $P$, without loss of generality, $P$ starts at a vertex of degree $3$, say $Pv_1v_2v_3$; see $G_{20},G_{21}$ of Fig. 3.4.
Without loss of generality, the edge $v_1v_2$ is negative and $v_2v_3$ is positive.
Now contracting the path $P$ into a single vertex say $u$, the resulting graph $\bar{B}$ will have rank $3$.
By Theorem \ref{rank-3}, $\bar{B}$ is a complete tripartite graph, that is $^u\bar{B}=\Theta(2,2,1)$, where the two triangles of $\bar{B}$ have the same balanceness.
By Theorem \ref{balance},  there exists a sign function $\eta$ defined on $V(\bar{B})$, such that all edges of $\bar{B}^\eta$ are positive or all edges of $\bar{B}^\eta$ are positive except the common edge of the two triangles.
Returning to the original graph, $^uB=G_{20}$ or $G_{21}$ of Fig. 3.4.
Define a sign function $\theta$ on $B$ such that $\theta(v_1)=\theta(v_2)=\theta(v_3)=\eta(u)$, and $\theta(v)=\eta(v)$ for any other vertex $v$.
Then $B^\theta$ is the graph $G_{20}$ or $G_{21}$ with edge labeling as in Fig. 3.4.\hfill$\blacksquare$

By Theorem \ref{main3}, we easily get a characterization of simple bicyclic graphs (or balanced bicyclic signed graphs) of order $n$ with nullity $n-5$.

\begin{coro} \label{main4}
Let $B$ be a simple bicyclic graph of order $n$.
Then $\eta(B)= n-5$  if and only if $B$ is one of the following graphs: $G_{11},G_{12}$ of Fig. 3.2, $G_{17},G_{18}, G_{19}$ of Fig. 3.4.
\end{coro}

\begin{center}

\includegraphics[scale=0.5]{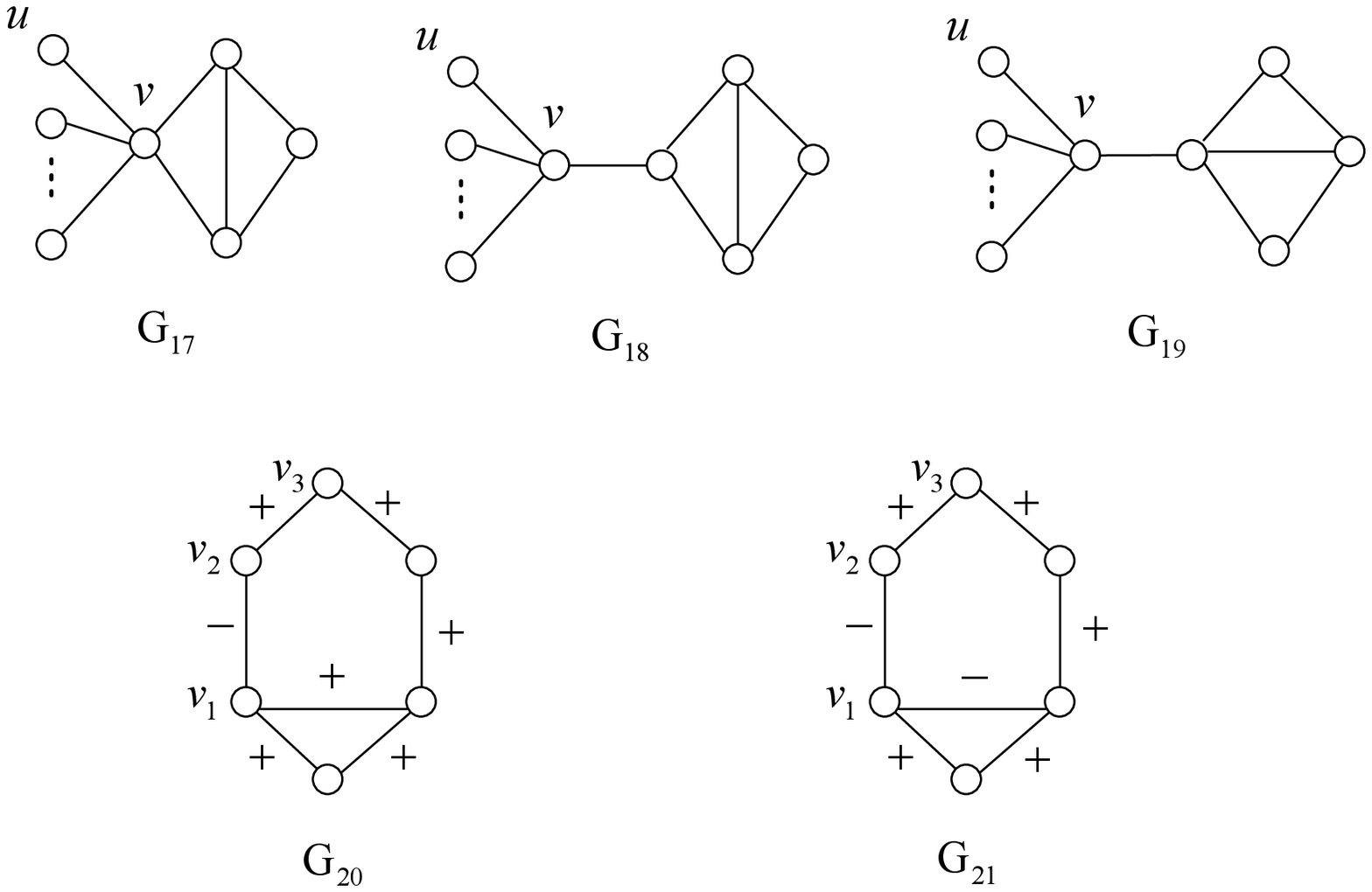}

{\small Fig. 3.4. Five  bicyclic graphs}
\end{center}

{\small

\end{document}